\documentclass[12pt]{amsart}
\usepackage{fullpage,url,amssymb}
\usepackage[all]{xy} 

\DeclareFontEncoding{OT2}{}{} 

\usepackage{color}

\newcommand{\GOOD}{\text{\sc Good}}
\newcommand{\Good}{\operatorname{Good}}
\newcommand{\Pfister}[1]{\langle \langle #1 \rangle \rangle}
\newcommand{\defi}[1]{\textsf{#1}}

\newcommand{\Aff}{{\mathbb A}}

\newcommand{\F}{{\mathbb F}}

\newcommand{\bbH}{{\mathbb H}}
\newcommand{\N}{{\mathbb N}}
\newcommand{\PP}{{\mathbb P}}
\newcommand{\Q}{{\mathbb Q}}

\newcommand{\Z}{{\mathbb Z}}

\newcommand{\kbar}{{\overline{k}}}

\newcommand{\calF}{{\mathcal F}}

\newcommand{\calT}{{\mathcal T}}

\newcommand{\OO}{{\mathcal O}}

\DeclareMathOperator{\newth}{th}
\DeclareMathOperator{\trdeg}{tr deg}

\DeclareMathOperator{\Spec}{Spec}
\DeclareMathOperator{\rk}{rk}
\DeclareMathOperator{\Char}{char}

\DeclareMathOperator{\End}{End}

\DeclareMathOperator{\cd}{cd}
\DeclareMathOperator{\vcd}{vcd}

\newcommand{\tors}{{\operatorname{tors}}}

\newcommand{\isom}{\simeq}

\newcommand{\Intersection}{\bigcap} 
\newcommand{\Union}{\bigcup} 
\newcommand{\union}{\cup} 
\newcommand{\tensor}{\otimes}

\newtheorem{theorem}{Theorem}[section]
\newtheorem{lemma}[theorem]{Lemma}
\newtheorem{corollary}[theorem]{Corollary}
\newtheorem{proposition}[theorem]{Proposition}

\theoremstyle{definition}
\newtheorem{definition}[theorem]{Definition}

\newtheorem{example/fact}[theorem]{Example/Fact}
\theoremstyle{remark}
\newtheorem{remark}[theorem]{Remark}

\usepackage[alphabetic,backrefs,lite]{amsrefs} 

\begin{document}

\title[Function fields over anti-Mordellic fields]{First-order definitions 
in function fields over anti-Mordellic fields}
\author{Bjorn Poonen}
\thanks{B.P. was supported by NSF grant DMS-0301280
        a Packard Fellowship, and the Miller Institute for Basic Research
        in Science.  He thanks the 
        Isaac Newton Institute for hosting a visit in the summer of 2005.}
\address{Department of Mathematics, University of California, 
        Berkeley, CA 94720-3840, USA}
\email{poonen@math.berkeley.edu}
\urladdr{http://math.berkeley.edu/\~{}poonen}
\author{Florian Pop}
\address{Department of Mathematics, University of Pennsylvania, DRL, 
        209 S 33rd Street, Philadelphia, PA 19104. USA}
\email{pop@math.upenn.edu}
\urladdr{http://math.penn.edu/\~{}pop}
\date{February 20, 2006}


\maketitle

\section{Introduction}\label{S:introduction}

\begin{definition}
A field $k$ is \defi{anti-Mordellic} (or \defi{large}) 
if every smooth curve with a $k$-point
has infinitely many $k$-points.
\end{definition}

Separably closed fields, henselian fields,
and real closed fields are examples of anti-Mordellic fields.
An algebraic extension of an anti-Mordellic field is anti-Mordellic: 
see \cite{Pop1996}, Proposition~1.2.

\begin{definition}
Let $k$ be a field.
A \defi{function field} over $k$ is a finitely generated extension $K$ of $k$
with $\trdeg(K|k)>0$.
\end{definition}

\begin{definition}
The \defi{constant field} of a field $K$ finitely generated over $k$ 
is the relative algebraic closure of $k$ in $K$.
\end{definition}

\begin{theorem}
\label{T:geometric constants}
There exists a formula $\phi(t)$
that when interpreted in a field $K$ finitely generated over
an anti-Mordellic field $k$ defines the constant field.
\end{theorem}

\begin{theorem}
\label{T:distinguishing classes}
For each of the following classes of fields,
there is a sentence that is true for fields in that class
and false for fields in the other five classes:
\begin{enumerate}
\item finite and anti-Mordellic fields
\item number fields
\item function fields over finite fields
\item function fields over anti-Mordellic fields of characteristic~$>0$
\item function fields over anti-Mordellic fields of characteristic~$0$
\item function fields over number fields
\end{enumerate}
\end{theorem}

\begin{remark}
It is impossible to distinguish finite fields from
anti-Mordellic fields with a single sentence,
since a nontrivial ultraproduct of finite fields
is anti-Mordellic.
\end{remark}

Finally, we have a few theorems characterizing algebraic dependence.
Some of these require that the ground field $k$ be 
``$2$-cohomologically well behaved'' in the sense of
Definition~\ref{D:well behaved} in Section~\ref{S:independence}.
The following theorems will be proved in Section~\ref{S:independence}.

\begin{theorem}
\label{T:alg dep}
There exists a formula $\phi_n(t_1,\ldots,t_n)$
such that for every $K$ finitely generated over a real closed
or separably closed field $k$,
and every $t_1,\ldots,t_n \in K$,
the formula holds if and only if $t_1,\ldots,t_n$ are algebraically
dependent over $k$.
\end{theorem}

\begin{theorem}
\label{T:alg dep over any well behaved field}
Let $k$ be a $2$-cohomologically well behaved field.
Let $K|k$ be a finitely generated extension.
Then there exists a first order formula (depending on $K$ and $k$)
with $r$ free variables,
in the language of fields augmented by a predicate for a subfield,
that when interpreted for elements $t_1,\ldots,t_r \in K$
with the subfield being $k$ 
holds if and only if the elements are algebraically independent over $k$.
\end{theorem}

\begin{corollary}
\label{C:alg dep over special fields}
Let $k$ be a finite field, a number field, or a
$2$-cohomologically well behaved anti-Mordellic field.
Then there exists a first order formula (depending on $K$ and $k$)
with $r$ free variables, in the language of fields,
that when interpreted for elements $t_1,\ldots,t_r \in K$
holds if and only if the elements are algebraically independent over $k$.
\end{corollary}

\begin{proof}
Theorem~1.4 of \cite{Poonen2005-uniform-preprint}
handles the case where $k$ is finite or a number field.
If $k$ is anti-Mordellic, 
combine Theorems \ref{T:geometric constants} 
and~\ref{T:alg dep over any well behaved field}.
\end{proof}

\begin{remark}
We do not know if Theorem~\ref{T:alg dep over any well behaved field}
and Corollary~\ref{C:alg dep over special fields}
can be made uniform in $k$ and $K$,
i.e., whether the formula can be chosen independent of $k$ and $K$.
\end{remark}

\section{Defining the constants}

In this section we prove Theorem~\ref{T:geometric constants}.

\begin{lemma}
\label{L:complement}
Let $k$ be an infinite field of characteristic $p$.
Let $S_0$ be a finite subset of $k$,
and let $S = \{s^{p^n}: s \in S_0, n \in \N\}$.
Then $k-S$ is infinite.
\end{lemma}

\begin{proof}
If $k$ is algebraic over $\F_p$, then $S$ is finite, so $k-S$ is infinite.
Otherwise, choose $t \in k$ transcendental over $\F_p$;
then for a given $s \in k$, the set $\{s^{p^n}: n \in \N\}$
contains at most one element of $\{t^\ell: \text{ $\ell$ is prime}\}$,
so $k-S$ is infinite.
\end{proof}

\begin{lemma}
\label{L:non-isotrivial}
Let $k$ be an infinite field.
Let $V$ be a $k$-variety.
Let $\{ C_a \}$ be a non-isotrivial family of curves of genus $\ge 2$ 
over $k$ with parameter $a$.
Then there exist infinitely many $a \in k$
such that all rational maps from $V$ to $C_a$ are constant.
\end{lemma}

\begin{proof}
Let $p$ be the characteristic of $k$.
A theorem of Severi \cite{Samuel1966}*{Th\'eor\`eme~2}
states that there are only finitely many fields $L$ between $k$
and the function field $K$ of $V$
such that $L$ is the function field of a curve of genus $\ge 2$ over $k$
and $K$ is separable over $L$.
Thus the set $S$ of $a \in k$ such that 
$C_a$ admits a non-constant rational map from $V$
is a finite set $S_0$ together with (if $p>0$) 
the $p^n$-th powers of the elements of $S_0$ for all $n \in \N$.
By Lemma~\ref{L:complement}, $k-S$ is infinite.
\end{proof}

\begin{proof}[Proof of Theorem~\ref{T:geometric constants}]
Without loss of generality we may assume that $k$ is relatively
algebraically closed in $K$.
The discriminant of $x^5+ax+1$ (with respect to $x$)
is $256 a^5+3125$; if $\Char k \notin \{2,5\}$, 
this is a nonconstant squarefree polynomial in $a$,
so the family of affine curves $C_a\colon y^2=x^5+a x+1$ 
has both smooth and nodal curves, and is therefore non-isotrivial.
If $\Char k=5$, the family $C_a\colon y^2=x^7+ax+1$ is non-isotrivial
for the same reason; and if $\Char k=2$, the family
$C_a\colon y^2+y=x^5+ax$ is non-isotrivial,
since a direct calculation (using the fact that the unique Weierstrass point
must be preserved) shows that no two members of this family
are isomorphic over an algebraic closure of $k$.
The projection $x \colon C_a \to \Aff^1$
is \'etale above $0 \in \Aff^1(k)$.

For $a \in K$, define
\[
  S_a:=\left\{\frac{x_1}{x_2} : (x_1,y_1), \, (x_2,y_2) \in C_a(K) 
       \text{ with $x_2 \ne 0$} \right\}.
\]
(A very similar definition was used in the proof 
of \cite{Koenigsmann2002}*{Theorem~2}.)
We have
\begin{enumerate}
\item\label{I:fact1}
If $a \in k$, then $k \subseteq S_a$.
{\em Proof:} Let $f(x,y)=0$ is the equation of $C_a$ in $\Aff^2$.
Let $c \in k$.
The map $(x_1,x_2)\colon C_a \times C_a \to \Aff^2$ is \'etale
above $(0,0)$, so the point 
$(x_1,y_1,x_2,y_2)=(0,1,0,1)$ 
on the inverse image $Y$
of the line $x_1=cx_2$ in $C_a \times C_a$ 
is smooth.
Since $k$ is anti-Mordellic, $Y$ has infinitely many
other $k$-points, so $c \in S_a$.
\item\label{I:fact2}
There exists $a_0 \in k$ such that $S_{a_0}=k$.
{\em Proof:} 
Let $V$ be an integral $k$-variety with function field $K$.
Lemma~\ref{L:non-isotrivial} gives $a_0 \in k$ such that
there is no nonconstant rational map $V \dashrightarrow C_{a_0}$ over $k$.
Equivalently, $C_a(K)=C_a(k)$.
So $S_{a_0} \subseteq k$, and we already know the opposite inclusion.
\item\label{I:fact3}
If $a \in K-k$, then $S_a$ is finite.
{\em Proof:} 
By the function field
analogue of the Mordell conjecture~\cite{Samuel1966}*{Th\'eor\`eme~4},
$C_a(K)$ is finite,
so $S_a(K)$ is finite.
\end{enumerate}

Let $A$ be the set of $a \in K$ such that $S_a$ is a field containing $a$.
Let $L := \Intersection_{a \in A} S_a$.
Then $L$ is uniformly definable by a formula.
By \eqref{I:fact3}, $A \subseteq k$ (a finite field cannot contain
an element transcendental over $k$).
Now by \eqref{I:fact1} and \eqref{I:fact2}, $L=k$.
\end{proof}

\begin{remark}
\label{R:S_a'}
Suppose $K$ is finitely generated over a field $k$,
and $k$ is relatively algebraically closed in $K$.
By the Weil conjectures applied to $Y$, 
there exists an explicit positive integer $m$
such that \eqref{I:fact1} is true also in the case where $k$ is a finite
field of size $>m$.
Let $S_a'$ be the union of $S_a$ with the set of zeros
of $x^q-x$ in $K$ for all $q \in \{2,3,\ldots,m\}$.
Let (1)', (2)', (3)' be the statements analogous to 
\eqref{I:fact1}, \eqref{I:fact2}, \eqref{I:fact3} 
but with $S_a'$ in place of $S_a$.
Then (1)', (2)', (3)' remain true 
for anti-Mordellic $k$, but now (1)' and (3)'
hold also for finite $k$.
\end{remark}

\section{Some facts about quadratic forms}
\label{S:facts}

\begin{proposition}
\label{P:Springer}
Let $q(x_1,\ldots,x_n)$ be a quadratic form over a field $K$,
and let $L$ be a finite extension of $K$ of odd degree.
If $q$ has a nontrivial zero over $L$,
then $q$ has a nontrivial zero over $K$.
\end{proposition}

\begin{proof}
This is well known: see \cite{LangAlgebra}*{Chapter~V, Exercise~28}.
\end{proof}

\begin{corollary}
\label{C:purely inseparable}
Let $K$ be a field of characteristic not $2$.
Let $q$ be a quadratic form over $K$.
Let $L$ be a purely inseparable extension of $K$.
If $q$ has a nontrivial zero over $L$,
then $q$ has a nontrivial zero over $K$.
\end{corollary}

\begin{proof}
If $q$ has a nontrivial zero over $L$,
the coordinates of this zero generate a finite purely inseparable extension
of $K$, so we may assume $[L:K]<\infty$.
Now the result follows from Proposition~\ref{P:Springer}.
\end{proof}

For nonzero $a$,
let $\Pfister{a}$ denote the quadratic form $x^2+ay^2$ 
and let 
$\Pfister{a_1,\ldots,a_n}=\Pfister{a_1}\tensor \cdots \tensor \Pfister{a_n}$ 
be the $n$-fold Pfister form.

\begin{lemma}
\label{L:local parameters}
Let $k$ be a field,
and let $V$ be an integral $k$-variety with function field $K$.
Suppose that $v$ is a regular point on $V$,
and that $t_1,\ldots,t_m$ are part of a system of local parameters at $v$.
Let $q$ be a diagonal quadratic form over $k$ 
having no nontrivial zero over the residue field of $v$.
Then $q \tensor \Pfister{t_1,\ldots,t_m}_d$ 
has no nontrivial zero over $K$.
\end{lemma}

\begin{proof}
This result is essentially contained in \cite{Pop2002}.
The proof is given again in Lemma~A.5 in \cite{Poonen2005-uniform-preprint}.
\end{proof}

\begin{lemma}
\label{L:small trdeg}
Let $\ell$ be a field of characteristic not $2$.
Let $L$ be a finitely generated extension of $\ell$.
Suppose that every $3$-fold Pfister form $\Pfister{a,b,c}$ over $L$
has a nontrivial zero.
Then
\begin{enumerate}
\item
$\trdeg(L|\ell) \le 2$.
\item
If moreover $L$ admits a valuation that is trivial on $\ell^\times$
such that $\ell$ maps isomorphically to the residue field,
and not every element of $\ell$ is a square in $\ell$,
then $\trdeg(L|\ell) \le 1$.
\end{enumerate}
\end{lemma}

\begin{proof}\hfill
(1) 
Let $t_1,\ldots,t_d$ be a transcendence basis for $L|\ell$.
Let $K$ be the maximal separable extension of $\ell(t_1,\ldots,t_d)$
contained in $L$.
Let $V$ be an integral variety over $\ell$ with function field $K$.
Replacing $V$ by an open subset if necessary,
we may assume that $(t_1,\ldots,t_d)\colon V \to \Aff^n_\ell$
is \'etale.
If $\ell$ is infinite, choose $(a_1,\ldots,a_d) \in \Aff^n(\ell)$ 
in the image of $V$; then by Lemma~\ref{L:local parameters},
$\Pfister{t_1-a_1,\ldots,t_d-a_d}$ has no nontrivial zero over $K$,
and hence by Corollary~\ref{C:purely inseparable}, 
no nontrivial zero over $L$.
If $\ell$ is finite, choose $(a_1,\ldots,a_d) \in \Aff^n(\ell')$
in the image of $V$ for some $\ell'|\ell$ of odd degree,
and repeat the previous argument with the minimal polynomial
$P_{a_i}(t_i)$ of $a_i$ over $\ell$ in place of $t_i-a_i$.
In either case, this Pfister form contradicts the hypothesis if $d \ge 3$.
Thus $d \le 2$.

(2)
Suppose not.
Then by (1), $\trdeg(L|\ell)=2$. By the resolution of singularities 
for surfaces (see e.g.\ \cite{Abhyankar1969}), we may choose a regular 
projective surface $V$ over $\ell$ with function field $L$.
The center of the given valuation on $V$ is an $\ell$-rational point 
$v \in V(\ell)$; hence $v$ is actually a smooth point of $V$. Choose 
local parameters $u_1,u_2$ at $v$. Let $\alpha \in \ell$ be a non-square.
By Lemma~\ref{L:local parameters}, $\Pfister{-\alpha,u_1,u_2}$ 
has no nontrivial zero over $L$. This contradicts the hypothesis.
\end{proof}

\begin{lemma}
\label{L:potential density}
Let $X$ be a variety over an infinite field $k$.
There exists an integer $m$
such that the points on $X$ of degree $\le m$ over $k$
are Zariski dense in $X$.
\end{lemma}

\begin{proof}
The desired property depends only on the birational class of $X$
over $\kbar$.
Therefore, enlarging $k$, we may reduce to the case where $X$ is a 
geometrically integral closed hypersurface in $\PP^n$.
Choose $P \in (\PP^n-X)(k)$.
Projection from $Q$ determines a generically
finite rational map from $X$ to $\PP^{n-1}$,
and the fibers above $k$-points in a Zariski dense open subset
of $\PP^{n-1}$ contain points of bounded degree.
These points are Zariski dense in $X$.
\end{proof}

\section{Distinguishing classes of fields}
\label{S:distinguishing}

\begin{proposition}
\label{P:finite and anti-Mordellic}
There is a sentence $\phi$ that is true for finite fields 
and anti-Mordellic fields,
false for function fields over any field,
and false for number fields.
\end{proposition}

\begin{proof}
Let $K$ be a field.
Define $S_a'$ as in Remark~\ref{R:S_a'}.
Let $\phi$ be the sentence saying that $S_a'=K$ for all $a \in K$.
This is true if $K$ is finite or anti-Mordellic.

If $K$ is a function field,
then (3)' (whose proof is valid over any $k$)
shows that for some $a$, the set $S_a'$ is finite.
If $K$ is a number field,
then $S_a'$ is finite for all but finitely many $a$,
by the Mordell conjecture~\cite{Faltings1983} applied to $C_a$.
In both these cases, there exists $a \in K$ with $S_a' \ne K$.
\end{proof}

We can generalize Theorem~\ref{T:geometric constants} 
to include finitely generated extensions of finite fields:

\begin{proposition}
\label{P:constants over finite and anti-Mordellic}
There exists a formula that for $K$ finitely generated
over a finite or anti-Mordellic field $k$
defines the constant field.
\end{proposition}

\begin{proof}
We may assume that $k$ is relatively algebraically closed in $K$.
We use the notation of 
the proof of Theorem~\ref{T:geometric constants} 
and Remark~\ref{R:S_a'}.
Let $A'$ be the set of $a \in K$ such that $S_a'$ is a field containing $a$.
Let $k_1:=\Intersection_{a \in A'} S_a'$.
Theorem~1.3 of \cite{Poonen2005-uniform-preprint} 
gives a formula that defines the constant subfield
if $K$ is finitely generated over a finite field;
over any field $K$, let $k_2$ be the subset it defines.
Define
\[
        \tilde{k} := \begin{cases}
                k_1, &\text{ if $S_a' \supseteq k_1$ for every $a \in k_1$,}\\
                k_2, &\text{ otherwise.}
        \end{cases}
\]
The subset $\tilde{k}$ is definable by a uniform formula;
we claim that $\tilde{k}=k$.

If $k$ is anti-Mordellic, then by 
the proof of Theorem~\ref{T:geometric constants},
$k_1=k$, and $\tilde{k}=k_1=k$.

Now suppose $k$ is finite, so $k_2=k$.
The set $k_1$ is a field (since it is an intersection of fields),
and it contains $k$ by Remark~\ref{R:S_a'}.
If $k_1=k$, then $\tilde{k}=k$.
If $k_1 \supsetneq k$, and $a \in k_1-k$,
then by \eqref{I:fact3}, 
$S_a'$ is finite, so it cannot contain $k_1$;
thus $\tilde{k}=k_2=k$.
\end{proof}

\begin{proposition}
\label{P:function field over finite and anti-Mordellic}
There exists a sentence that is true for function fields 
over finite or anti-Mordellic fields
and false for number fields and function fields over number fields.
\end{proposition}

\begin{proof}
Use the sentence that says that the formula in 
Proposition~\ref{P:constants over finite and anti-Mordellic}
defines a field satisfying the sentence
of Proposition~\ref{P:finite and anti-Mordellic}.
\end{proof}

\begin{proposition}
\label{P:function field over finite field}
There is a sentence in the language of rings extended by a unary predicate
that when interpreted in a function field $K$ over a field $k$
(not necessarily relatively algebraically closed)
with the unary predicate defining $k$
is true if and only if $k$ is finite.
\end{proposition}

\begin{proof}
By \cite{Poonen2005-uniform-preprint}*{Remark~5.1},
there is a formula $\phi(x,y)$ in the language of rings
such that when it is interpreted in a function field $K$
with finite constant field $\ell$,
\[
        \{y \in K: \phi(x,y) \} = \ell[x]
\]
for each $x \in K$.
By \cite{Rumely1980}, there is a formula $\psi$ defining a family
of subsets
that when interpreted in $\ell(x)$ for $\ell$ finite
gives exactly the family of nontrivial valuation rings in $\ell(x)$
(possibly with repeats).

Now let $K$ be a function field over an arbitrary field $k$.
We claim that $k$ is finite if and only if for some $x \in K$ 
the following hold:
\begin{enumerate}
\item\label{I:condition-ring}
The set $R$ defined by $\phi(x,\cdot)$ is a ring containing $k$ and $x$.
\item\label{I:condition-valuationrings}
The family $\calF$ defined by $\psi$ interpreted over the fraction 
field $L$ of $R$ defines a set of nontrivial valuation rings in $L$,
each containing $k$.
\item\label{I:condition-constants}
The intersection of the valuation rings in $\calF$ is a field $\ell$.
\item\label{I:condition-x}
The element $x$ is not in $\ell$.
\item\label{I:condition-residuefield}
The field $\ell$ maps isomorphically to the residue field
of some valuation ring in $\calF$.
\item\label{I:condition-char2}
If $2=0$, then $[L:L^2]=2$.
\item\label{I:condition-Pfister}
If $2 \ne 0$, then every $3$-fold Pfister form $\Pfister{a,b,c}$ over $L$
has a nontrivial zero, and some element of $\ell$ is not a square in $\ell$.
\item\label{I:condition-intersection}
The intersection of the rings in $\calF$ containing $R$ equals $R$.
\item\label{I:condition-PID}
Every ideal $aR+bR$ of $R$ generated by two elements is principal.
\item\label{I:condition-irreducibles}
The elements $x-a$ for $a \in \ell$ are irreducible, 
and generate pairwise distinct ideals of $R$.
\item\label{I:condition-divisibility}
There exists a nonzero $f \in R$ 
divisible in $R$ by $x-a$ for all $a \in \ell$.
\end{enumerate}
(These conditions can be expressed by a first order sentence
in the language of rings with a predicate for $k$.)

If $k$ is finite, and $x \in K$ is not in the constant field $\ell$ of $K$,
then $R=\ell[x]$ for a finite field $\ell$, and 
conditions \eqref{I:condition-ring}--\eqref{I:condition-divisibility} hold.

Conversely, suppose that 
conditions \eqref{I:condition-ring}--\eqref{I:condition-divisibility} hold
for some $x \in K$.
If $\Char K=2$, then \eqref{I:condition-char2} implies $\trdeg(L|\ell) \le 1$.
If $\Char K\ne 2$, then 
\eqref{I:condition-residuefield}
and \eqref{I:condition-Pfister} imply that $\trdeg(L|\ell) \le 1$,
by Lemma~\ref{L:small trdeg}.
Thus in every case, $\trdeg(L|\ell) \le 1$.
By \eqref{I:condition-constants}, $\ell$ is an intersection of valuation
rings, so it is relatively algebraically closed in $L$.
By \eqref{I:condition-x}, $x \in L-\ell$, so $\trdeg(L|\ell) = 1$.
Since $L$ is a function field over $k$ and $k \subseteq \ell$,
$L$ is a function field of transcendence degree $1$ over $\ell$.
By \eqref{I:condition-intersection}, $R$ is integrally closed in $L$;
in particular it contains the integral closure $R_0$ of $\ell[x]$ in $L$.
Thus $R_0$ is a Dedekind domain with fraction field $L$.
Any ring between a Dedekind domain and its fraction field is itself
a Dedekind domain, so $R$ is a Dedekind domain.
By \eqref{I:condition-PID}, $R$ is a principal ideal domain,
and hence a unique factorization domain.
Now \eqref{I:condition-irreducibles} and \eqref{I:condition-divisibility} 
imply that $\ell$ is finite.
So $k$ is finite.
\end{proof}

\begin{proposition}
\label{P:finite versus anti-Mordellic}
There is a sentence that is true for function fields over finite fields
and false for function fields over anti-Mordellic fields.
\end{proposition}

\begin{proof}
Combine Propositions \ref{P:constants over finite and anti-Mordellic}, 
and~\ref{P:function field over finite field}.
\end{proof}

\begin{proposition}
\label{P:characteristic of function field over anti-Mordellic}
There exists a sentence that for a function field $K$ over
a finite or anti-Mordellic field
is true if and only if $\Char K=0$.
\end{proposition}

Before beginning the proof of 
Proposition~\ref{P:characteristic of function field over anti-Mordellic},
we need a few definitions and a lemma.
If $M$ is an Abelian group and $n \ge 1$,
let $M[n]$ be the kernel of the multiplication-by-$n$ map $M \to M$.
Also define $M_\tors:=\Union_{n \ge 1} M[n]$.
If $E\colon y^2=f(x)$ is an elliptic curve 
over a field $K$ of characteristic $\ne 2$,
and $t \in K$,
then the twisted elliptic curve $E_t$ is defined
by $f(t) y^2=f(x)$ over $K$.
We will use the following, 
which is essentially a special case of a result of Moret-Bailly.

\begin{lemma}
\label{L:Moret-Bailly}
Let $k$ be a field of characteristic~$0$.
Let $K$ be a function field over $k$.
Let $E \colon y^2=f(x)$ be an elliptic curve over $k$,
where $f$ is a cubic polynomial.
Then there are infinitely many $t \in K$
with $f(t) \in K^{\times} - k^\times K^{\times 2}$
such that $E_t(K)$ is a finitely generated Abelian group
with $\rk E_t(K) = \rk \End_K(E)$.
\end{lemma}

\begin{proof}
We may enlarge $k$ to assume that $K$ is the function field
of a geometrically irreducible curve over $k$.
Replacing $f(x)$ by $f(x+c)$ for suitable $c \in k$,
we may assume that $f(0) \ne 0$.

We use the terminology in \cite{Moret-Bailly2005preprint}*{\S1.5}.
Let $\Gamma$ be the smooth projective model of
the $y^2=x^4 f(1/x)$; cf. \cite{Moret-Bailly2005preprint}*{1.4.5(ii)}.
By \cite{Moret-Bailly2005preprint}*{2.3.1},
there exists $g \in K-k$ that is admissible for $\Gamma$.
By \cite{Moret-Bailly2005preprint}*{1.8(ii) and~1.4.7}, 
$\GOOD(k) \cap \Z$ is infinite.

We claim that for any $\lambda \in \GOOD(k) \cap \Z$,
the value $t:=\frac{1}{\lambda g}$ satisfies the required conditions.
For such $\lambda$ and $t$,
we have $\lambda \in \Good(k)$ 
by \cite{Moret-Bailly2005preprint}*{1.5.4(i)};
thus $E'\colon (\lambda g)^4 f(\frac{1}{\lambda g}) y^2 = f(x)$
is an elliptic curve over $K$
such that $E'(K)$ is finitely generated 
and $\rk E'(K) = \rk \End_K(E)$.
By definition, $E'$ is isomorphic to $E_t$.

Let $K \kbar$ be a compositum of $K$ with an algebraic closure of $\kbar$
over $k$.
If $f(t)$ were in $k^\times K^{\times 2}$,
then $E'$ would be isomorphic over $K \kbar$ to $E$,
so $E'(K\kbar) \isom E(K \kbar) \supseteq E(\kbar)$ 
would not be finitely generated,
contradicting the definition of $\GOOD(k)$.
\end{proof}

\begin{proof}[Proof of 
Proposition~\ref{P:characteristic of function field over anti-Mordellic}]
Use $\neg \phi$, where $\phi$ is a sentence
equivalent to the following:
$2=0$ or there exists an extension $L$ of $K$ with $[L:K] \le 2$
such that for $\ell$ the subset defined in by the formula
of Proposition~\ref{P:constants over finite and anti-Mordellic}
applied to $L$,
there exist distinct $e_1,e_2,e_3 \in \ell$
such that if we write $f(x):=(x-e_1)(x-e_2)(x-e_3)$,
then for all $t \in L$ 
with $f(t) \in L^\times - \ell^\times L^{\times 2}$,
the twist $E_t$ of $E\colon y^2=f(x)$
satisfies $\# E_t(L)/2E_t(L) \ge 64$.
For the $K$ we are interested in,
$L$ is a function field over a finite or anti-Mordellic field,
so $\ell$ is the constant field of $L$.

If $\Char K=2$, then $\phi$ is true.
Now suppose $K$ is a function field over an anti-Mordellic field
of characteristic $p>2$.
Let $L$ be a compositum of $K$ with $\F_{p^2}$.
Let $E$ be an elliptic curve over $\F_p$ with $\#E(\F_p)=p+1$.
Then the $p^2$-Frobenius endomorphism of $E$ is 
multiplication by $-p$,
so $\rk \End_{\F_{p^2}}(E)=4$, and $E[2] \subseteq E(\F_{p^2})$.
The curve $E_{\F_{p^2}}$ has an equation $y^2=f(x)$
where $f(x):=(x-e_1)(x-e_2)(x-e_3)$ 
with distinct $e_1,e_2,e_3 \in \F_{p^2} \subseteq \ell$.
Suppose $t \in L$ satisfies $f(t) \in L^\times - \ell^{\times} L^{\times 2}$.
The restriction on $t$ implies that $E_t$ is not isomorphic over $L$
to an elliptic curve over $\ell$, so $E_t(L)$ is finitely generated.
Quadratic twists of an elliptic curve have the same endomorphism ring,
so the ring $\OO:=\End_L(E_t)$ is a maximal order in a non-split
quaternion algebra $\bbH$ over $\Q$.
Since $E_t(L) \tensor \Q$ is an $\bbH$-vector space,
$4 \mid \rk_\Z E_t(L)$.
The point $(t,1) \in E_t(L)$ has infinite order,
since under the $L(\sqrt{f(t)})$-isomorphism $E_t \to E$ 
mapping $(x,y)$ to $(x,y\sqrt{f(t)})$
it corresponds to a point of $E$
whose $x$-coordinate is transcendental over $\ell$.
Thus $\rk_\Z E_t(L)>0$, so $\rk_\Z E_t(L) \ge 4$.
Also, $E_t[2] \subseteq E_t(L)$, 
so $\# E_t(L)/2E_t(L) \ge 2^2 \cdot 2^4 = 64$.

Now suppose that $K$ is a function field over an anti-Mordellic field
of characteristic~$0$.
Suppose $L$ is an extension with $[L:K] \le 2$,
and $e_1,e_2,e_3 \in \ell$ are distinct.
By Lemma~\ref{L:Moret-Bailly} applied to $L$ over $\ell$,
there exists $t \in L$ with $f(t) \notin \ell^\times L^{\times 2}$
such that $E_t(L)$ is finitely generated with
$\rk E_t(L) = \rk \End_L(E)$.
Since $\rk \End_L(E) \in \{1,2\}$, and since $E_t(L)_\tors$
is generated by at most $2$ elements,
we get $\# E_t(L)/2E_t(L) \le 2^2 \cdot 2^2 = 16$.
\end{proof}

\begin{proof}[Proof of Theorem~\ref{T:distinguishing classes}]
Taking $d=0$ in the first claim of Theorem~1.5(3) of \cite{Pop2002}
gives a sentence that is true for number fields
and false for function fields over number fields.
Combining this with
Propositions \ref{P:finite and anti-Mordellic},
\ref{P:function field over finite and anti-Mordellic},
\ref{P:finite versus anti-Mordellic},
and \ref{P:characteristic of function field over anti-Mordellic}
gives the result.
\end{proof}

\section{Detecting algebraic dependence}
\label{S:independence}

We begin by recalling the following general facts: Let $E$ be 
an arbitrary field of characteristic $\neq2$. In particular, 
$\mu_2=\{\pm1\}$ is contained in $E$. We denote by $G_E$ the 
absolute Galois group of $E$, and view $\mu_2$ as a $G_E$-module. 

1) Let $\cd^0_2(E) \in \N \union \{\infty\}$ 
be the supremum over all the natural 
numbers $n$ such that ${\rm H}^n(E,\mu_2)\neq 0$. Since the 
$2$-cohomological dimension $\cd_2(E)$ is defined similarly, but
the supremum is taken over all possible $2$-torsion $G_k$-modules,
one has  
\[
   \cd^0_2(E)\leq\cd_2(E)\,. 
\]
Also define $\vcd_2(E):=\cd_2(E(\sqrt{-1}))$.
   
2) Recall the Milnor Conjecture (proved by Voevodsky et al.) 
It asserts that the $n^{\newth}$ cohomological invariant 
$e_n \colon I_n(E)/I_{n+1}(E)\to{\rm H}^n(E,\mu_2)$, which maps each 
$n$-fold Pfister form $\Pfister{a_1,\ldots,a_n}$ to the cup product 
$(-a_1)\cup\dots\cup (-a_n)$, is a well defined isomorphism. Using 
the Milnor Conjecture one can describe $\cd^0_2(E)$ via the 
behavior of Pfister forms as follows: $n>\cd^0_2(E)$ if and 
only if every $n$-fold Pfister form over $E$ 
represents $0$ over $E$.

3) There exists a field $E$ with $\cd^0_2(E)<\cd_2(E)$.
For instance, let $E$ be a maximal pro-$2$ Galois extension of a
global or local field of characteristic $\neq2$. Then every 
element of $E$ is a square, so $\cd^0_2(E)=0$. On the other 
hand, since the Sylow $2$-groups of $G_E$ are non-trivial, one 
has $\cd_2(E)>0$
by \cite{SerreGaloisCohomology}*{\S I.3.3, Corollary~2}.

\begin{definition} 
\label{D:well behaved}
A field $E$ is said to be \defi{$2$-cohomologically well behaved}
if $\Char E \ne 2$ and
for every finite extension $E'|E$ containing $\sqrt{-1}$ 
one has $\cd^0_2(E')=\cd_2(E')<\infty$.
\end{definition}

\begin{remark}
\label{R:finite extension}
If $E$ is $2$-cohomologically well behaved, then
and $E'|E$ is a finite extension containing $\sqrt{-1}$,
then
\[
        \cd^0_2(E')=\vcd_2(E')=\vcd_2(E),
\]
since $\cd_2(E')=\cd_2(E(\sqrt{-1}))$ 
by \cite{SerreGaloisCohomology}*{\S II.4.2, Proposition~10}.
\end{remark}

\begin{example/fact}
The following fields, when of characteristic $\neq 2$, 
are $2$-cohomologically well behaved:
\begin{itemize}
\item
separably closed fields (trivial),
\item
finite fields (follows from \cite{SerreGaloisCohomology}*{II.\S3}),
\item
local fields (follows from \cite{SerreGaloisCohomology}*{II.\S4.3}),
\item
number fields (follows from \cite{SerreGaloisCohomology}*{II.\S4.4}), and
\item
finitely generated fields 
(follows from the above and Proposition~\ref{cd20 of function field} below).
\end{itemize}
\end{example/fact}

\begin{proposition}
\label{cd20 of function field}
If $E$ is $2$-cohomologically well behaved, 
and $E'$ is a function field over $E$,
then $E'$ is $2$-cohomologically well behaved
and $\vcd_2(E')=\vcd_2(E)+\trdeg(E'|E)$.
\end{proposition}

\begin{proof}
We may assume $\sqrt{-1} \in E$.
The case $\trdeg(E'|E)=0$ follows from Remark~\ref{R:finite extension}.
By induction on $\trdeg(E'|E)$, it will suffice to prove
that $\cd_2^0(E')=\vcd_2(E)+1$ for every extension $E'|E$
with $\trdeg(E'|E)=1$.
We may assume that $E'$ is separably generated over $E$.
Let $X$ be a curve over $E$ with function field $E'$, 
let $P$ be a smooth point on $X$,
let $\kappa$ be the residue field of $P$,
and let $t \in E'$ be a uniformizer at $P$.
Let $n=\cd_2^0(\kappa)=\vcd_2(E)$.
By definition, there exists an $n$-fold Pfister form 
$\Pfister{\bar{a}_1,\ldots,\bar{a}_n}$ 
that does not represent $0$ over $\kappa$.
Lift each $\bar{a}_i$ to an $a_i$ in the local ring at $P$.
Then $\Pfister{a_1,\ldots,a_n,t}$ does not represent $0$ over $E'$.
Thus $\cd_2^0(E') \ge \vcd_2(E)+1$.
On the other hand, $\cd_2^0(E') \le \vcd_2(E') = \vcd_2(E)+1$
by \cite{SerreGaloisCohomology}*{\S II.4.2, Proposition~11},
so we have equality.
\end{proof}

\begin{proposition}
\label{P:more-general-alg-dep}
Let $k$ be a field which is $2$-cohomologically well behaved, 
and let $e=\vcd_2(k)$. Let $K|k$ be a finitely generated extension.
Then the following hold:
\begin{enumerate}
\item 
For each $n \in \Z_{\ge 0}$,
there exists a sentence $\phi_n$ in the 
language of fields (depending on~$e$) 
such that $\phi_n$ is true in $K$ if and only if $\trdeg(K|k)=n$.

One can take $\phi_n$ to be the following sentence: Every 
$(e+n+1)$-fold Pfister form over $K[\sqrt{-1}]$ represents $0$,
but there exist $(e+n)$-fold Pfister forms over $K[\sqrt{-1}]$ 
which do not represent $0$.

\item 
For elements $t_1,\ldots,t_r \in K^\times$,
the following are equivalent:
\begin{enumerate}
\item $(t_1,\ldots,t_r)$ are algebraically independent over $k$. 
\item 
There exists a finite separable extension $l|k$ 
(depending on $t_1,\ldots,t_r$) containing $\sqrt{-1}$
and elements $a_1,\dots,a_e,b_1,\dots,b_r \in l^\times$ such 
that $\Pfister{a_1,\dots,a_e,t_1-b_1,\dots t_r-b_r}$ 
does not represent $0$ over $Kl$.
\end{enumerate}
\end{enumerate} 
\end{proposition}

\begin{proof}\hfill

(1) By the discussion preceding Proposition~\ref{P:more-general-alg-dep}
we have: 
\[
  \cd^0_2(K[\sqrt{-1}])=\cd^0_2(k[\sqrt{-1}])+\trdeg(K|k)=e+\trdeg(K|k).
\] 
Now use the characterization of $\cd^0_2$ in terms of Pfister forms.

(2), (b)~$\Rightarrow$~(a): Suppose for the sake of obtaining
a contradiction that $(t_1,\dots,t_r)$ is algebraically dependent 
over $k$. Let $L=l(t_1,\dots,t_r)\subset Kl$. Since $\sqrt{-1}\in l$, 
we have: 
\[
   \cd_2(L)=\cd_2(l)+\trdeg(L|l)=e+\trdeg(L|l)<e+d
\]   
Thus by the discussion above, every $(e+d)$-fold Pfister form over
$L$ represents $0$ over $L$. In particular, for all $(a_i)_i$ 
and $(b_j)_j$ as in (b), the resulting $(e+d)$-fold Pfister form
$\Pfister{a_1,\dots,a_e,t_1-b_1,\dots t_r-b_r}$ represents $0$ over 
$L$. Since $L\subseteq Kl$, it follows that 
$\Pfister{a_1,\dots,a_e,t_1-b_1,\dots t_r-b_r}$ represents 
$0$ over $Kl$, a contradiction!

(2), (a) $\Rightarrow$ (b): The proof is an adaptation 
from and similar to \cite{Pop2002}, Section~1. 
By extending the list $\calT:=(t_1,\ldots,t_r)$,
we may assume that it is a transcendence basis for $K|k$.
Let $K_0|k(\calT)$ be the 
relative separable closure of $k(\calT)$ in $K$. Thus $\calT$ 
is a separable transcendence basis of $K_0|k$, and $K|K_0$ is a 
finite purely inseparable field extension. 
Further let $R$ be the integral closure of $k[\calT]$ in $K_0$,
and let $X=\Spec R$.
The $k$-embedding 
$k[\calT]\hookrightarrow R$ defines a finite $k$-morphism
$\phi \colon X \to \Spec k[\calT]=\Aff^r_k$. Further, since 
$K_0|k(\calT)$ is separable, the $k$-morphism $\phi$
is generically \'etale. Therefore, $\phi$ is \'etale on a Zariski 
dense open subset $U\subset X'$. We choose a finite separable
extension $l|k$ containing $\sqrt{-1}$
such that $U(l)$ is non-empty. 
Choose $x\in U(l)$, and let $b:=(b_1,\dots,b_r)=\phi(x)$ 
be its image in $\Aff^r_k(l)=l^r$. 
Then $t_1-b_1,\dots,t_r-b_r$ are local parameters at $x$.
Since $\cd_2^0 l = e$, we may choose $a_1,\ldots,a_e \in l^\times$
such that $\Pfister{a_1,\dots,a_e}$ has no nontrivial zero over $l$.
Then by Lemma~\ref{L:local parameters}, 
$\Pfister{a_1,\ldots,a_e,t_1-b_1,\ldots, t_r-b_r}$ 
has no nontrivial zero over $K_0 l$.
By Corollary~\ref{C:purely inseparable},
$\Pfister{a_1,\ldots,a_e,t_1-b_1,\ldots, t_r-b_r}$ 
has no nontrivial zero over $K l$.
\end{proof}

\begin{proof}[Proof of Theorem~\ref{T:alg dep over any well behaved field}] 
Theorem~1.4 of \cite{Poonen2005-uniform-preprint}
handles the case where $k$ is finite, so assume that $k$ is infinite.
By replacing $k$ with a finite extension $k'$ and simultaneously
replacing $K$ with $Kk'$ (these extensions can be interpreted 
over $(K,k)$), we may assume that $K$ is the function field of a 
geometrically integral variety $X$ over $k$ where $\sqrt{-1} \in k$, 
and by Lemma~\ref{L:potential density} we may assume that the points 
of degree $\le m$ on $X$ are Zariski dense. Now, by the same proof 
as in Proposition~\ref{P:more-general-alg-dep}(2), $t_1,\ldots,t_r$ 
are algebraically independent over $k$ if and only if there exists 
an extension $l|k$ of degree $\le m$ such that there exist 
$a_1,\ldots,a_e,b_1,\ldots,b_r \in l^\times$ such that 
$\Pfister{a_1,\ldots,a_e,t_1-b_1,\ldots, t_r-b_r}$ has no 
nontrivial zero over $Kl$. The preceding statement is expressible 
as a certain first order formula evaluated at $t_1,\ldots,t_r$.
\end{proof}

Unfortunately, in the case $\Char=2$ we do not have at our 
disposal an easy way to relate $\trdeg(K|k)$ to (some) well 
understood invariants (say similar to the cohomological
dimension). In the case $k$ is separably closed, one can though 
employ the theory of $C^{(p)}_i$ fields. Recall that a field $E$ 
is said to be a $C^{(p)}_i$ field, if every system of homogeneous 
forms 
\[f_\rho(X_1,\dots,X_n)\quad(\rho=1,\dots,r)
\] 
has a non-trivial common zero, provided the degrees $d_\rho$ of 
the forms satisfy: $n>\sum_\rho d^i_\rho$ and $(p,d_\rho)=1$ for 
all $\rho$.

The following are well known facts about $C^{(p)}_i$ fields, see
e.g., \cite{Pfister1995}:

1) Suppose that $E$ is a $p$-field, i.e., 
every finite extension $E'|E$ has degree a power of $p$. 
Then $E$ is a $C^{(p)}_0$ field.

2) If $E$ is a $C^{(p)}_i$ field, then every finite extension 
$E'|E$ is again a  $C^{(p)}_i$ field.

3) If $E$ is a $C^{(p)}_i$ field, then the rational function 
field $E(t)$ in one variable over $E$ is an $C^{(p)}_{i+1}$ 
field.

In particular, if $k$ is a $C^{(p)}_i$ field, and $K|k$ is a 
function field with $\trdeg(K|k)=d$, then $K$ is a $C^{(p)}_{i+d}$ 
field.  

Now let $K|k$ be a function field. For every $\ell>0$ and every 
system $\underline t=(t_1,\dots,t_r)$ of elements of $K^\times$, 
let
\[
 q^{(\ell)}_{(t_1,\dots,t_r)}=
    \sum_{\underline i}\underline t^{\underline i}X^\ell_{\underline i}
 \] 
be the ``generalized Pfister form" of degree $\ell$ in $\ell^r$ 
variables as introduced in \cite{Pop2002}, Section~1, p.~388. Here 
${\underline i}$ is a multi-index ${\underline i}=(i_1,\dots,i_r)$, 
with $0\leq i_j<\ell$. 
 
\begin{proposition}
\label{P:less-gen-than-more-general}
Let $k$ be a $p$-field. Let $K|k$ be a function field. 
For every $\ell\neq\Char(K)$ that is relatively prime to $p$, 
and every system $(t_1,\dots,t_r)$ 
of elements of $K^\times$, let $q^{(\ell)}_{(t_1,\dots,t_r)}$ 
be the corresponding form of degree $\ell$ in $\ell^r$ variables 
over $K$. Then one has:
\begin{enumerate}
\item For every $r>\trdeg(K|k)$, and every $(t_1,\dots,t_r)$, the 
resulting form $q^{(\ell)}_{(t_1,\dots,t_r)}$ represents $0$ 
over $K$. 
\item For a given system $(t_1,\dots,t_r)$ the following conditions 
are equivalent: 
\begin{enumerate}
\item $(t_1,\dots,t_r)$ is algebraically independent over $k$.

\item there exist $b_1,\dots,b_r \in k$ such that
$q^{(\ell)}_{(t_1-b_1,\dots,t_r-b_r)}$ does not represent $0$
over $K$.
\end{enumerate}

\item 
In particular, for each $n \in \Z_{\ge 0}$
there exists a sentence in the language of fields 
that holds for $K$ if and only if $\trdeg(K|k)=n$. 
\end{enumerate}  

Thus given algebraically independent elements $x_1,\dots, x_r\in K$ 
over $k$, the relative algebraic closure $L$ of $k(x_1,\dots,x_r)$ 
in $K$ is described by a predicate in one variable $x$ as follows:
\[
  L=\{\,x\in K\mid (x_1,\dots,x_r,x)
       \hbox{ is not algebraically independent over $k$}\,\}
\]
\end{proposition}

\begin{proof}[Proof of 
              Proposition~\ref{P:less-gen-than-more-general}] \hfill

(1): By the discussion above, $K$ is a $C^{(p)}_d$ field for $d=\trdeg(K|k)$. 

(2), (b)~$\Rightarrow$~(a): Let $L=k(t_1,\dots,t_r)$. 
If $t_1,\ldots,t_r$ are algebraically dependent,
then $\trdeg(L|k)<r$, so by (1), any form
$q^{(\ell)}_{(t_1-b_1,\dots,t_r-b_r)}$ represents $0$ over $L$,
and hence represents $0$ over $K$.

(2), (a) $\Rightarrow$ (b): The proof is very similar to
the proof of the corresponding implication in 
Proposition~\ref{P:more-general-alg-dep}. 
The relative separable closure $K_0$ of $k(t_1,\ldots,t_r)$ in $K$
is the function field of an \'etale cover $U \to \Aff^r_k$
with the morphism being given by $(t_1,\ldots,t_r)$.
Choose $(b_1,\ldots,b_n) \in \Aff^r(k)$ and a closed point $u \in U$
above it.
If $l$ is the residue field of $u$, then $K_0$ embeds into
the iterated Laurent power series field 
$\Lambda:=l((t_r-b_r))\dots((t_1-b_1))$,
and $K$ embeds into a purely inseparable 
finite extension $\Lambda'$ of $\Lambda$.
The field $\Lambda$ has a natural valuation $v$ whose value group
is $\Z^r$ ordered lexicographically, generated 
by $v(t_i-b_i)$ for $1\leq i\leq r$.
The values of the coefficients of 
$q^{(\ell)}_{(t_1-b_1,\dots,t_r-b_r)}$ 
are distinct modulo $\ell$
(they even a system of representatives for $\Z^r/\ell \Z^r$).
If we extend $v$ to a valuation $v'$ on $\Lambda'$,
the value group $G$ of $v'$ contains $\Z^r$ with index prime to $\ell$,
so the $v'$-valuations of these coefficients
have distinct images in $G/\ell G$.
Now for any non-zero system
of elements $\underline x=(x_{\underline i})_{\underline i}$ from 
$K \subseteq \Lambda'$,  
$q^{(\ell)}_{(t_1-b_1,\dots,t_r-b_r)}(\underline x)$ 
is a sum of elements having distinct $v'$-valuations
(distinct even modulo $\ell$).
So $q^{(\ell)}_{(t_1-b_1,\dots,t_r-b_r)}(\underline x)\neq 0$.
 
The remaining assertions of Proposition~\ref{P:less-gen-than-more-general}
are clear. 
\end{proof}

\begin{proof}[Proof of Theorem~\ref{T:alg dep}] \hfill

\underline{Case 1}: $\Char(k)\neq2$.

If $k$ is either real closed or separably closed, then $l:=k[\sqrt{-1}]$
is the unique finite separable field extension of $k$ containing $\sqrt{-1}$.
Thus the result follows from Proposition~\ref{P:more-general-alg-dep}~(2).

\underline{Case 2}: $\Char(k)=2$.

Then $k$ is a $2$-field, so it is a $C{(2)}_0$ field. 
To conclude, one applies Proposition~\ref{P:less-gen-than-more-general}
with $p=2$ and $\ell=3$. 
\end{proof}

\section*{Acknowledgments} 

We thank Laurent Moret-Bailly for some discussions 
of his paper~\cite{Moret-Bailly2005preprint}.

\end{document}